\documentclass{article}

\usepackage{arxiv}
\usepackage{amssymb}
\usepackage{amsmath}
\usepackage{empheq}
\usepackage{amscd}
\usepackage{graphicx}
\usepackage{graphics}
\usepackage[noadjust]{cite}
\usepackage{amsthm}
\usepackage{caption}
\usepackage{multirow}
\usepackage{array}
\usepackage{rotating}

\usepackage{indentfirst}
\usepackage{tikz}
\usepackage{calc}
\usepackage{epsfig}
\usepackage{natbib}
\usepackage[makeroom]{cancel}
\usepackage{multicol}
\usepackage{csquotes}
\usepackage{enumitem}

\usepackage{hyperref}       

\begin{document}
\title{Extended full waveform inversion in the time domain by the augmented Lagrangian method}

\author{\href{https://orcid.org/0000-0002-9879-2944}{\includegraphics[scale=0.06]{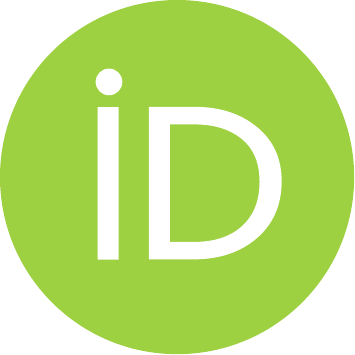}\hspace{1mm}Ali Gholami} \\
  Institute of Geophysics, University of Tehran, Tehran, Iran.
  \texttt{agholami@ut.ac.ir} \\ 
   \And
 \href{http://orcid.org/0000-0003-1805-1132}{\includegraphics[scale=0.06]{orcid.pdf}\hspace{1mm}Hossein S. Aghamiry} \\
  University Cote d'Azur - CNRS - IRD - OCA, Geoazur, Valbonne, France. 
  \texttt{aghamiry@geoazur.unice.fr} 
  \And
\href{http://orcid.org/0000-0002-4981-4967}{\includegraphics[scale=0.06]{orcid.pdf}\hspace{1mm}St\'ephane Operto} \\ 
  University Cote d'Azur - CNRS - IRD - OCA, Geoazur, Valbonne, France. 
  \texttt{operto@geoazur.unice.fr}
  }


\renewcommand{\shorttitle}{Data Reconstruction Inversion ~~~~~~~~~~~~~~~~~~~~~~~~~~~~~~~~~~~~~~~~~~~~~~ Gholami et al.}

\maketitle
\begin{abstract}
Extended full-waveform inversion (FWI) has shown promising results for accurate estimation of subsurface parameters when the initial models are not sufficiently accurate. Frequency-domain applications have shown that the augmented Lagrangian (AL) method
solves the inverse problem accurately with a minimal effect of the penalty parameter choice. 
Applying this method in the time domain, however, is limited by two main factors:
(1) The challenge of data-assimilated wavefield reconstruction due to the lack of an explicit time stepping and
(2) The need to store the Lagrange multipliers, which is not feasible for the field-scale problems. 
We show that these wavefields are efficiently determined from the associated data (projection of the wavefields onto the receivers space) by using explicit time stepping. Accordingly,  based on the augmented Lagrangian, a new algorithm is proposed which performs in ``data space" (a lower dimensional subspace of the full space) in which the wavefield reconstruction step is replaced by reconstruction of the associated data, thus requiring optimization in a lower dimensional space (convenient for handling the Lagrange multipliers). 
We show that this new algorithm can be implemented efficiently in the time domain with existing solvers for the FWI and at a cost comparable to that of the FWI while benefiting from the robustness of the extended FWI formulation.
The results obtained by numerical examples show high-performance of the proposed method for large scale time-domain FWI. 
\end{abstract}
\graphicspath{{"./figures/"}}
\section{Introduction}
Full-waveform inversion (FWI) is the state-of-the-art data-fitting based imaging method for subsurface parameters estimation. 
The parameters of interest appear as coefficients in system of hyperbolic partial differential equations (PDE) which serves as the forward modeling engine for simulating seismic waves. 
Sampling these waves at receiver locations gives us synthetic seismograms. 
FWI tries to find the set of parameters from which we can compute synthetic seismograms which best fit (in the least-squares sense) the
observed seismograms (data) \citep{Pica_1990_NIS,Virieux_2009_OFW}.
The forward problem is nonlinear thus iterative linearization methods are natural choice to solve the problem \citep{Tarantola_2005_IPT}.

Traditional formulation of FWI eliminates the intermediate variable (wavefields) from the equations and represents observed data directly in terms of the Green's functions which are nonlinearly related to the model parameters.
The standard gradient (steepest descent) method solves this nonlinear inverse problem iteratively. 
At each iteration, the parameters are updated along the steepest-descent direction.
This direction is efficiently determined by using adjoint-state method \citep{Plessix_2006_RAS} that is a simple zero-lag cross-correlation between two calculated wavefields: one is the forward wavefield in the current model and the other is a (backward) wavefield obtained by back-propagation of the data residuals acting as if they were sources at receiver locations.
A more accurate search direction is obtained by using Newton or quasi-Newton methods by including curvature information via the Hessian matrix or an approximation of it \citep{Metivier_2017_TRU}. A simple while effective approximation of the Hessian matrix can also be obtained by  using only its diagonal terms that are calculated by zero-lag correlation of the forward wavefield \citep{Shin_2001_IAP}. 
Calculated step direction is then scaled by an appropriate step length that can be chosen using a line search method, such as the backtracking line search, and added to the previous model to get the next iterate. 
This process is repeated until the convergence is achieved. 
In this way, the quality of the search direction is the main factor governing the algorithm convergence and it is a directly influenced by the quality of forward and backward wavefields. 

In contrast to the FWI, extended FWI brings flexibilities to fit the data by softening the wave equation constraint at least at early iterations \citep{vandenBerg_2001_CSI,VanLeeuwen_2013_MLM,Huang_2018_SEW,Aghamiry_2019_IWR}. 
The associated objective function is bivariate and besides model parameters we need to also manipulate data-assimilated wavefields \citep{Aghamiry_2020_AED}, specific wavefields that are obtained by solving the data equation and the wave equation simultaneously.
The objective function is biconvex and the augmented Lagrangian (AL) method provides an excellent framework for its minimization with a minimal sensitivity to the penalty parameter choice \citep{Gabay_1976_ADA,Nocedal_2006_NO,Aghamiry_2019_IWR}.
However, there are two main obstacles that limit the usage of this method in the time domain \citep{Wang_2017_RFI,Rizzuti_2019_ADF,vanLeeuwen_2019_ANO,Aghamiry_2020_AED}.
The first is reconstruction  of data-assimilated wavefields and the second is the need to store the Lagrange multipliers, which are of the size of the wavefields and thus infeasible for the field-scale problems. 

In this paper, we provide a new formulation for solving the AL-based extended FWI. 
An equivalent formulation of the iteratively-refined wavefield reconstruction inversion  \citep{Aghamiry_2019_IWR} is derived in the data space which require manipulating Lagrange multipliers of the data size. Furthermore, it replaces the wavefield reconstruction problem in the full space with a residual data reconstruction problem in the data space.
The resulting reconstruction system is of the data size and thus more tractable compared with the wavefield reconstruction in the original formulation. 
Based on the fact that the AL algorithm is convergent for a wide range of penalty parameters \citep{Gabay_1976_ADA,Nocedal_2006_NO,Aghamiry_2019_IWR}, we derive a simplified algorithm which only requires four wave equation solves at each iteration which can be carried out by usual time-stepping. 
An interesting observation is that the AL-based extended FWI updates the model parameters via a normalized adjoint-state formulation similar to the quasi-Newton algorithm for the classical FWI. 
The only difference is in the forward and backward wavefields which are cross-correlated. 
This can be thought of as a generalized algorithm which can be applied to both the extended and reduced FWI and in both the time and frequency domains.

The contents of this paper are organized as follows: We begin with the general form of the FWI problem as a partial-differential equation (PDE) constrained data-fitting problem.
We then review  the two approaches of classical FWI and the extended FWI for solving this optimization problem. 
A brief review of the augmented Lagrangian (AL) method is presented with its advantages and limitations for time domain FWI.
Then we derive an algorithm for AL-based extended FWI and show how it can be efficiently applied in the time domain. 
Finally, we will demonstrate the efficiency of the proposed algorithm by using  Marmousi II and BP 2004 benchmark velocity models.

\section{Theory}
FWI seeks the subsurface parameters by solving the following PDE constrained optimization problem \citep{VanLeeuwen_2013_MLM,Aghamiry_2019_IWR}:
\begin{align} \label{main_optim}
& \min_{\bold{m},\bold{u}} ~~\frac12\|\bold{Pu-d}\|_2^2~~~ \text{subject to}~ ~~\bold{A(m)u=b},
\end{align} 
where $\|\cdot\|_2$ denotes the Euclidean norm, $\bold{u}\in {\mathbb{R}}^{N\times 1}$ is the wavefield,
$\bold{d} \in {\mathbb{R}}^{M \times 1}$ is the recorded seismic data generated by the sampling operator $\bold{P} \in {\mathbb{R}}^{M\times N}$ (with $M\ll N$), $\bold{A(m)} \in {\mathbb{R}}^{N \times N}$ is the PDE operator, and $\bold{b} \in {\mathbb{R}}^{N\times 1}$ is the source term. 
More specifically, $N=N_x\times N_z \times N_t$ and $M=N_t\times N_r$ with $N_x$ and $N_z$ the number of spatial grid samples of the model parameter $\bold{m}$ (square of slowness) in horizontal and vertical dimensions, $N_t$ the number of time samples, and $N_r$ is the number of receivers.
%
\subsection{FWI}
Eliminating the wavefield from equation \ref{main_optim} leads to the conventional (reduced) formulation of FWI as \citep{Pratt_1998_GNF}
\begin{align} \label{main_RDS}
\min_{\bold{m}} ~~\frac12\|\bold{G(m)}\bold{b}-\bold{d}\|_2^2.
\end{align}
where $\bold{G(m)}=\bold{P}\bold{A(m)}^{-1}$ denotes the forward modeling operator.  
This objective function is very hard to solve due to the presence of the inverse operator $\bold{A(m)}^{-1}$ (which includes the Green's functions as its columns). 
The standard gradient-based iterative methods solve this nonlinear inverse problem iteratively \citep{Tarantola_1988_TBI,Pica_1990_NIS,Pratt_1998_GNF,Virieux_2009_OFW}, in which at iteration $k$th the parameters are updated as
\begin{equation}
\bold{m}_{k+1} = \bold{m}_k + \alpha_k\delta\bold{m}_k,
\end{equation}
where $\alpha_k$ is a step length that can be chosen using a line search method, such as the backtracking line search, and $\delta\bold{m}_k$ is a step direction. This step direction is selected as opposite gradient of the objective function ($\bold{g}_k$, in the steepest descent method) or determined by solving a Newton or quasi-Newton system:
\begin{equation} \label{Newton}
\delta\bold{m}_k = -\bold{H}_k^{-1}\bold{g}_k,
\end{equation}
where $\bold{H}_k$ is the Hessian matrix or an approximation of it. Both the gradient and the Hessian are evaluated at the current model $\bold{m}_k$. The gradient can be computed efficiently by using the adjoint state method \citep{Pica_1990_NIS,Plessix_2006_RAS}.
\begin{equation} \label{FWR_grad}
\bold{g}_k = \int_0^{\tau} [\ddot{\bold{u}}^r_k(t)\circ \bold{v}^r_k(t)]   dt,
\end{equation}
where  $\ddot{\bold{u}}^r_k(t)$ is the second time derivative of the forward wavefield induced by the source and $\bold{v}^r_k(t)$ is the adjoint wavefield obtained by back-propagation of the data residuals $\delta\bold{d}^r_k$ (which acts as a source emitting in a time-reversed manner).
\begin{equation} \label{delta_d}
\delta\bold{d}^r_k=\bold{d}-\bold{G}(\bold{m}_k)\bold{b}.
\end{equation}
In equation \ref{FWR_grad}, both wavefields are propagated in the current model $\bold{m}_k$ and they are recorded over the time interval $(0, \tau)$. Also, $\circ$ denotes the element-wise product operator. 

A simple while effective approximation of the Hessian or the pseudo-Hessian includes only the diagonal terms 
 \citep{Shin_2001_IAP}, which gives
\begin{equation}  \label{FWI_adjoint_state}
\delta\bold{m}_k = -\frac{\int_0^{\tau} [\ddot{\bold{u}}^r_k(t)\circ \bold{v}^r_k(t)]  dt}{\int_0^{\tau} |\ddot{\bold{u}}^r_k(t)|^2 dt},
\end{equation}
where the division and the square-power are element-wise operations.
This indeed is the zero-lag crosscorrelation image between the forward and adjoint wavefields (imaging condition) normalized by the square of the source illumination strength \citep{Claerbout_1971_TUT}.
A more accurate Hessian matrix can also be implement to improve the model update by solving equation \ref{Newton} via an iterative algorithm while computing Hessian-vector products by second-order adjoint state formulas \citep{Metivier_2017_TRU}.

However, FWI requires us to provided an initial model of sufficient quality for convergence to the desired solution \citep{Virieux_2009_OFW}. The next section provides a remedy for this issue using extended FWI.
In the formulas above, superscript $r$ shows that the variable corresponds to the reduce FWI. We will use superscript $e$ to refer to the extended FWI.
\subsection{Extended FWI}
Extended FWI methods have been developed to decrease the nonlinearity of the objective function and hence increase the robustness of local optimization methods to the initial model \citep{vandenBerg_2001_CSI,VanLeeuwen_2013_MLM,Huang_2018_SEW,Aghamiry_2019_IWR}. 
The simple form of the extended formulation can be obtained by the penalty formulation of the original equation \ref{main_optim}:
\begin{equation} \label{primal_penalty}
\min_{\bold{m,u}^e} ~\frac{1}{2}\|\bold{Pu}^e-\bold{d}\|_2^2 + \frac{\mu}{2}\|\bold{A(m)}\bold{u}^e-\bold{b}\|_2^2,
\end{equation}
where $\mu>0$ is the penalty parameter. 
This extended formulation can be thought of as an ``error-in-variable" method which accounts for errors in the wave equation, and it reduces to the FWI objective in equation \ref{main_RDS} as $\mu \rightarrow \infty$.
Equation \ref{primal_penalty} is biconvex and it is usually solved by alternating minimization \citep{Beck_2017_FOM}.
This algorithm however suffers from two major drawbacks which limit its application in a time domain practical setting.
The first is manipulating the u-subproblem (computation and storage of the wavefield for a given model) and the other is the issue related to the penalty parameter choice.

The issue of penalty parameter can be overcome by adding a Lagrangian term of the form $\boldsymbol{\lambda}^T(\bold{A(m)}\bold{u}^e-\bold{b})$ to the penalty function \citep{Nocedal_2006_NO,Aghamiry_2019_IWR}:
\begin{equation} \label{primal_AL}
\min_{\bold{m,u}^e}\max_{\boldsymbol{\lambda}} ~\text{Penalty}(\bold{m,u}^e) +\boldsymbol{\lambda}^T(\bold{A(m)}\bold{u}^e-\bold{b}).
\end{equation}
where superscript $T$ denotes adjoint operator, Penalty denotes the penalty objective function in equation \ref{primal_penalty} and $\boldsymbol{\lambda}$ is the Lagrange multiplier for the PDE constraint.
The alternating direction method of multipliers \citep{Gabay_1976_ADA,Goldstein_2009_SBM,Boyd_2011_DOS,Aghamiry_2019_IWR} solves the augmented Lagrangian in equation \ref{primal_AL} iteratively by choosing $\bold{u}$ and $\boldsymbol{\lambda}$ and then minimizing the objective as a function of $\bold{m}$. The resulting $\bold{m}$ is used to choose a new $\bold{u}$ and $\boldsymbol{\lambda}$, and the process repeats until the convergent is achieved.
\begin{subequations} \label{main_RI_WRI}
\begin{align} 
\bold{m}_{k}&= \arg\min_{\bold{m}} ~\frac{1}{2}\|\bold{A}(\bold{m})\bold{u}^e_{k-1}-\bold{b} - \boldsymbol{\lambda}_{k-1}\|_2^2, \label{main_RI_WRI_m}\\
\bold{u}^e_{k}&= \arg\min_{\bold{u}} ~\frac{1}{2}\|\bold{Pu}-\bold{d}\|_2^2 + \frac{\mu}{2}\|\bold{A}_{k}\bold{u}-\bold{b}-\boldsymbol{\lambda}_{k-1}\|_2^2, \label{main_RI_WRI_u}\\
\boldsymbol{\lambda}_{k} &= \boldsymbol{\lambda}_{k-1} + \bold{b} - \bold{A}_{k}\bold{u}^e_{k},  \label{main_RI_WRI_d}
\end{align}
\end{subequations}
beginning with initial $\bold{u}_0$ (which is built from an initial model) and $\boldsymbol{\lambda}_0=\bold{0}$,  where $\bold{A}_{k}\equiv \bold{A}(\bold{m}_{k})$ and $k$ is the iteration count. 
Note that \citet{Aghamiry_2019_IWR} first optimize in the wavefield $\bold{u}^e$ and then in the model parameters $\bold{m}$  but for the method to develop in this paper the permuted order in equation \ref{main_RI_WRI} is more convenient. This does not affect the final solution because the order of the optimization in biconvex problems can be permuted. 
The iteration in equation \ref{main_RI_WRI} without the ``add-back-the residual" term reduces to the wavefield reconstruction inversion (WRI) \citep{VanLeeuwen_2013_MLM}. This term accounts for the Lagrangian and besides increasing the convergence rate, reduces the sensitivity of the algorithm to the penalty parameter choice by ``adding-back the source residual".
The advantage of the augmented Lagrangian method is thus to use a small $\mu$ (to make the optimization well posed) and ``add-back-the-residual" \citep{Goldstein_2009_SBM}.
\citet{Aghamiry_2019_IWR} consider this iteration while solving the feasibility form of equation \ref{main_optim} in which the data are also updated by adding back the data residual. 
\begin{subequations} \label{main_RI_WRI_bd}
\begin{align} 
\bold{m}_{k}&= \arg\min_{\bold{m}} ~\frac{1}{2}\|\bold{A}(\bold{m})\bold{u}^e_{k-1}-\bold{b}_{k-1}\|_2^2, \label{main_RI_WRI_bd_m}\\
\bold{u}^e_{k}&= \arg\min_{\bold{u}} ~\frac{1}{2}\|\bold{Pu}-\bold{d}_{k-1}\|_2^2 + \frac{\mu}{2}\|\bold{A}_{k}\bold{u}-\bold{b}_{k-1}\|_2^2, \label{main_RI_WRI_bd_u}\\
\bold{d}_{k} &= \bold{d}_{k-1} + \bold{d} - \bold{P}\bold{u}^e_{k},  \label{main_RI_WRI_bd_d}\\
\bold{b}_{k} &= \bold{b}_{k-1} + \bold{b} - \bold{A}_{k}\bold{u}^e_{k},  \label{main_RI_WRI_bd_b}
\end{align}
\end{subequations}
beginning with  $\bold{b}_0=\bold{b}$ and $\bold{d}_0=\bold{d}$ which have the role of scaled Lagrange multipliers.
This algorithm has been successfully applied in the frequency-domain FWI where the most challenging subproblem, equation \ref{main_RI_WRI_u}, has a size of $N_x\times N_z$ and thus can be solved directly by sparse solvers. The Lagrangian multiplier of size $N_x\times N_z$ can also be conveniently stored during the iterations.
In the time domain, however, solving equation \ref{main_RI_WRI_u} is more challenging because the size of the system would be prohibitively large ($N_x\times N_z\times N_t$) due to the extra temporal dimension. Furthermore, $\bold{b}_{k}$ of the size $N_x\times N_z\times N_t$ would be dense; which causes serious storage problems. 
Several attempts have been make by different authors to solve these issues \citep{Wang_2017_RFI,Rizzuti_2019_ADF,vanLeeuwen_2019_ANO,Aghamiry_2020_AED}.
In the next section, we propose an algorithm for efficient treatment of these issues which performs in the data space instead of the full wavefield space.

\section{Data Reconstruction Inversion}
The general form of data-assimilated wavefield reconstruction, equation \ref{main_RI_WRI_u}, is 
\begin{equation} \label{primal_obj}
\min_{\bold{u}^e} ~\frac{1}{2}\|\bold{Pu}^e-\bold{d}\|_2^2 + \frac{\mu}{2}\|\bold{Au}^e-\bold{b}\|_2^2.
\end{equation}
Minimization of equation \ref{primal_obj} gives $\bold{u}^e$ as the solution of the following system:
\begin{equation} \label{Ex_u}
(\bold{P}^T\bold{P} + \mu\bold{A}^T\bold{A})\bold{u}^e=\bold{P}^T\bold{d} + \mu\bold{A}^T\bold{b},
\end{equation}
This system can be re-expressed by decoupling $\bold{Au}^e$ from other terms, making it suitable for time-stepping:
\begin{equation} \label{twosideeq}
\bold{Au}^e=\bold{b} + \frac{1}{\mu} \bold{G}^T\bold{(d-Pu}^e), 
\end{equation}
where $\bold{G}^T\equiv \bold{A}^{-T}\bold{P}^T= (\bold{A}^{-1})^T\bold{P}^T$.
However, since the unknown wavefield $\bold{u}^e$ appears in both sides of this equation, researchers tried different approaches to solve it. For example, \citet{Wang_2016_FIR} replaced the wavefield at the right-hand-side with the wavefield obtained by the reduced FWI. Later, \citet{Aghamiry_2020_AED} developed an algorithm which solves equation \ref{twosideeq} iteratively, beginning with the wavefield constructed by \citet{Wang_2016_FIR}.

In this paper, we solve equation \ref{twosideeq} in a two step procedure. Noting that $\bold{u}^e$  appears in the right-hand side only as (WRI predicted data) $\bold{Pu}^e$ thus if we have an estimate of $\bold{Pu}^e$ then $\bold{u}^e$  can be determined by solving equation \ref{twosideeq}. 
Accordingly, we rewrite this equation as
\begin{align} \label{twostepeq}       
\bold{Au}^e &=\bold{b} + \frac{1}{\mu} \bold{G}^T\delta\bold{d}^e, \\
          &=\bold{b} + \delta\bold{b}^e, \nonumber
\end{align}
where
\begin{equation} \label{residuals}
\delta\bold{b}^e=\bold{Au}^e-\bold{b}~~\text{and}~~\delta\bold{d}^e=\bold{d-Pu}^e,
\end{equation}
are respectively the source residual and data residual (of WRI). Note that the source residual is zero for the reduced FWI.
Up to now we have reduced the optimization dimension from $N$ (directly searching for $\bold{u}^e$ in equation \ref{Ex_u} or \ref{twosideeq}) to $M$ (searching for the residual vector $\delta\bold{d}^e$ in equation \ref{twostepeq}). 
However, the question is now ``what is the optimization problem for $\delta\bold{d}^e$?"

By eliminating $\bold{u}^e$ from equation \ref{residuals} we get the following relation between the two residual vectors $\delta\bold{b}^e$ and $\delta\bold{d}^e$:
\begin{equation} \label{residuals2}
\bold{G}\delta\bold{b}^e+\delta\bold{d}^e=\delta\bold{d}^r.
\end{equation}
where $\delta\bold{d}^r$  is the data residual of the FWI, equation \ref{delta_d}.

Furthermore, comparing equations \ref{twosideeq} and \ref{residuals} shows that
\begin{equation} \label{explicit_x}  
\delta\bold{b}^e=\frac{1}{\mu}\bold{G}^T\delta\bold{d}^e.
\end{equation}
Substituting $\delta\bold{b}^e$ from equation \ref{explicit_x} into equation \ref{residuals2} immediately gives an explicit $M\times M$ system for $\delta\bold{d}^e$ in the data space: 
\begin{equation} \label{explicit_y}
\left(\bold{G}\bold{G}^{T} + \mu\bold{I}\right)\delta\bold{d}^e=\mu\delta\bold{d}^r.
\end{equation}
This equation clearly shows that the FWI data residual, $\delta\bold{d}^r$, is a blurred version of the WRI data residual, $\delta\bold{d}^e$, in which the blurring kernel consists of the correlation of receiver-side Green's functions.
The solution of equation \ref{explicit_y} is also the minimizer of the following quadratic objective function:
\begin{equation} \label{delta_d_obj}
\min_{\delta\bold{d}^e} ~\frac{1}{2}\|\delta\bold{d}^e-\delta\bold{d}^r\|_2^2 + \frac{1}{2\mu}\|\bold{G}^T\delta\bold{d}^e\|_2^2.
\end{equation}

Plugging $\delta\bold{b}^e$ and $\delta\bold{d}^e$  from equations \ref{explicit_x} and \ref{explicit_y} into the primal objective in equation \ref{primal_obj} gives an equivalent objective for the extended FWI which is a weighted norm of the FWI data residual:
\begin{equation}
\frac{1}{2\mu}\|\delta\bold{d}^e\|_{\bold{Q}(\bold{m})}^2=\frac{\mu}{2}\|\delta\bold{d}^r\|_{\bold{Q}(\bold{m})^{-1}}^2.
\end{equation}
where $\|\bold{y}\|_{\bold{Q}(\bold{m})}^2=\bold{y}^T\bold{Q}(\bold{m})\bold{y}$ and $\bold{Q}(\bold{m})$ is the coefficient matrix in equation \ref{explicit_y}. 
This equivalence has also been confirmed differently by \citet{vanLeeuwen_2019_ANO} and \citet{Symes_2020_WRI}.

\subsection{IR-WRI Iteration in Data Space}
Using the data space formulation given above the IR-WRI iteration given in equation \ref{main_RI_WRI_bd} can be performed more efficiently in the data space. 
From equations \ref{main_RI_WRI_bd_d}, \ref{main_RI_WRI_bd_b} and \ref{residuals} we have that
\begin{align} \label{eqval}
\begin{cases}
\bold{d}_{k} =\bold{d}_{k-1} + \bold{d} - \bold{P}\bold{u}_{k} =\bold{d}+\delta\bold{d}^e_{k},\\
\bold{b}_{k} = \bold{b}_{k-1} + \bold{b} - \bold{A}_{k}\bold{u}_{k} =\bold{b} -\delta\bold{b}^e_{k}.
\end{cases}
\end{align} 
Furthermore, from equations  \ref{explicit_x} and \ref{explicit_y} we have
\begin{equation} \label{yk_1}
\begin{cases}
\delta\bold{b}^e_{k}=\frac{1}{\mu}\bold{G}_k^{T}\delta\bold{d}^e_{k},\\
\delta\bold{d}^e_{k}=\mu\bold{Q}_k^{-1}[\bold{d}_{k-1}-\bold{G}_k\bold{b}_{k-1}],
\end{cases}
\end{equation}
where
\begin{equation}
\bold{Q}_k=\bold{G}_k\bold{G}_k^{T}+\mu\bold{I}.
\end{equation}
Accordingly the wavefield $\bold{u}^e_k$ in equation \ref{main_RI_WRI_m} can be written as (see equation  \ref{twostepeq})
\begin{align} 
\bold{u}^e_{k} &=\bold{A}_k^{-1}[\bold{b}_{k-1}+ \delta\bold{b}^e_{k}] \\
 &=\bold{A}_k^{-1}[\bold{b} - \delta\bold{b}^e_{k-1}+ \delta\bold{b}^e_{k}] \\
&=\bold{A}_k^{-1}[\bold{b} +\frac{1}{\mu}\bold{G}_k^{T}(\delta\bold{d}^e_{k}-\delta\bold{d}^e_{k-1})],
\end{align}
(while replacing $\bold{G}_{k-1}$ with its update $\bold{G}_{k}$)
based on which the model update in equation \ref{main_RI_WRI_m} simplifies to
\begin{equation} \label{DRI_00}
\bold{m}_{k+1}= \arg\min_{\bold{m}} ~\frac{1}{2}\|\bold{A}(\bold{m})\bold{A}_{k}^{-1}\bold{b}_k^+-\bold{b}_k^-\|_2^2,
\end{equation}
in which
\begin{equation}  \label{DRI_11}
\begin{cases}
\bold{b}_k^+ =\bold{b} +\frac{1}{\mu}\bold{G}_k^{T}(\delta\bold{d}^e_{k}-\delta\bold{d}^e_{k-1}), \\
\bold{b}_k^-=\bold{b}- \frac{1}{\mu}\bold{G}_k^{T}\delta\bold{d}^e_k.
\end{cases}
\end{equation}
From equation \ref{delta_d_obj}, we see that $\delta\bold{d}_k^e$ is the minimizer of the following objective function
\begin{equation} \label{delta_d_obj2}
f(\delta\bold{d}^e)=\frac{1}{2}\|\delta\bold{d}^e-\delta\bold{d}^r_{k-1}\|_2^2 + \frac{1}{2\mu}\|\bold{G}_k^{T}\delta\bold{d}^e\|_2^2.
\end{equation}
Differentiating this objective with respect to $\delta\bold{d}^e$ gives 
\begin{equation} \label{delta_d_obj2_grad2}
\nabla  f(\delta\bold{d}^e)=\frac{1}{\mu}\bold{Q}_k\delta\bold{d}^e - \delta\bold{d}^r_{k-1} 
=\frac{1}{\mu}\bold{Q}_k\delta\bold{d}^e - \bold{d}_{k-1} + \bold{G}_k\bold{b}_{k-1}.
\end{equation}
Substituting the values of $\bold{d}_{k-1}$ and $\bold{b}_{k-1}$ from equations \ref{eqval} and \ref{yk_1} into equation \ref{delta_d_obj2_grad2} gives
\begin{align}
\nabla  f(\delta\bold{d}^e)
&=\frac{1}{\mu}\bold{Q}_k\delta\bold{d}^e - \bold{d}- \delta\bold{d}^e_{k-1}+\bold{G}_k(\bold{b} -\delta\bold{b}^e_{k-1}),\\
&=\frac{1}{\mu}\bold{Q}_k\delta\bold{d}^e-\delta\bold{d}_k^r-\delta\bold{d}^e_{k-1}-\bold{G}_k\delta\bold{b}^e_{k-1},\\
&=\frac{1}{\mu}\bold{Q}_k\delta\bold{d}^e -\delta\bold{d}_k^r-\delta\bold{d}^e_{k-1}-\frac{1}{\mu}\bold{G}_k\bold{G}_{k-1}^{T}\delta\bold{d}^e_{k-1},\label{needed}\\
&=\frac{1}{\mu}\bold{Q}_k\delta\bold{d}^e-\delta\bold{d}_k^r-\frac{1}{\mu}\bold{Q}_k\delta\bold{d}^e_{k-1}.\label{objgrad}
\end{align} 
Note that, in equation \ref{needed}, we replaced $\bold{G}_{k-1}$ with its update $\bold{G}_{k}$. 
Setting the gradient, equation \ref{objgrad}, equal to zero and solving for $\delta\bold{d}^e$ gives the exact minimizer, $\delta\bold{d}_k^e$, of the objective 
\begin{equation}
\delta\bold{d}_k^e=\delta\bold{d}^e_{k-1}+\mu\bold{Q}_k^{-1}\delta\bold{d}_k^r.\label{data_update}
\end{equation}
It is seen that in the augmented Lagrangian method the data residual vector is updated via a Newton like algorithm, equation \ref{data_update}, in which the search direction is defined by the FWI data residual $\delta\bold{d}_k^r$.
$\bold{Q}_k$ has the role of Hessian matrix. 

After introducing the auxiliary variable $\bold{y}_k=\delta\bold{d}^e_k/\mu$ and $\delta\bold{y}_k=\bold{y}_k-\bold{y}_{k-1}$, 
 equations \ref{DRI_00} and \ref{DRI_11} simplify to the following iteration

\begin{equation}  \label{main_PMS}
\begin{cases}
\delta\bold{y}_k=\bold{Q}_k^{-1}\delta\bold{d}^r_k,\\
\bold{y}_{k} =\bold{y}_{k-1}+\delta\bold{y}_k,\\
\bold{b}_k^+ =\bold{b} + \bold{G}_k^{T}\delta\bold{y}_k, \\
\bold{b}_k^-=\bold{b}- \bold{G}_k^{T}\bold{y}_{k},\\
\bold{m}_{k+1}= \arg\min_{\bold{m}} ~\frac{1}{2}\|\bold{A}(\bold{m})\bold{A}_{k}^{-1}\bold{b}_k^+-\bold{b}_k^-\|_2^2,
\end{cases}
\end{equation}
where $\bold{G}_k^{T}=\bold{A}_k^{-T} \bold{P}^T$. 
As seen, the wavefield reconstruction step in equation \ref{main_RI_WRI_bd} is replaced by a residual data reconstruction step, $\delta\bold{y}_k$ in equation \ref{main_PMS},  and thus we call this new iteration data reconstruction inversion (DRI). DRI only needs to store the dual variable $\bold{y}_k$ lying in the data space.

\subsection{DRI with Gradient Descent}
A fundamental property of all augmented Lagrangian based algorithms is that they are convergent for an approximate solution of the subproblems \citep{Gabay_1976_ADA,Boyd_2011_DOS,Yin_2013_EFG,Gholami_2015_NMI,Gholami_2017_CNA,Aghamiry_2019_IWR}. 
This motivated us to perform the minimization of objective $f$ in equation \ref{delta_d_obj2} approximately via a gradient step rather than full minimization given in equation \ref{data_update}.
\begin{equation} \label{data_update_grad}
\delta\bold{d}^e_{k} = \delta\bold{d}^e_{k-1}-\alpha_k\nabla f(\delta\bold{d}^e_{k-1}).
\end{equation} %
where $\alpha_{k}$ is a step length and $\nabla f$ is the gradient of the objective function determined in equation \ref{objgrad}, evaluated at $\delta\bold{d}^e_{k-1}$. 
From equation \ref{objgrad} we get that 
\begin{equation}
\nabla f(\delta\bold{d}^e_{k-1}) = -\delta\bold{d}_k^r.
\end{equation}
Accordingly, $\delta\bold{y}_k$ in equation \ref{main_PMS} becomes
\begin{equation} \label{hess_aprox}
\delta\bold{y}_k \approx \alpha_{k} \delta\bold{d}^r_k,
\end{equation}
which (after replacing $\bold{y}_{k}$ and $\bold{y}_{k-1}$ by $\bold{y}_{k}/\alpha_k$ and $\bold{y}_{k-1}/\alpha_k$) simplifies the DRI iteration as
\begin{equation}  \label{main_PMS2}
\begin{cases}
\bold{y}_{k} =\bold{y}_{k-1}+\delta\bold{d}^r_{k},\\
\bold{b}_k^+ =\bold{b} + \alpha_k\bold{G}_k^{T}\delta\bold{d}^r_k, \\
\bold{b}_k^-=\bold{b}- \alpha_k\bold{G}_k^{T}\bold{y}_{k},\\
\bold{m}_{k+1}= \arg\min_{\bold{m}} ~\frac{1}{2}\|\bold{A}(\bold{m})\bold{A}_{k}^{-1}\bold{b}_k^+-\bold{b}_k^-\|_2^2.
\end{cases}
\end{equation} 

\subsubsection{Step length determination}
In the DRI algorithm, equation \ref{main_PMS2}, the model parameters are updated by assuming a large penalty parameter and thus more weight given to minimization of the source residuals rather than the data residual. Accordingly, we determine the value of $\alpha_k$ by minimizing the data residuals:
\begin{equation}
\alpha_k=\arg\min_{\alpha}\|\bold{P}\bold{u}_k^e(\alpha)-\bold{d}\|_2^2,
\end{equation}
giving that
\begin{equation} \label{alpha}
\alpha_k=\frac{\bold{q}_k^T\delta\bold{d}^r_k}{\bold{q}_k^T\bold{q}_k},
\end{equation}
where $\bold{q}_k=\bold{G}_k\bold{G}_k^{T}\delta\bold{d}^r_k$.

\subsection{Update of the Model Parameters}
We consider constant density acoustic wave equation as the forward modeling operator, given by
\begin{equation} \label{PDE}
\bold{m}\circ \ddot{\bold{u}}(t)-\boldsymbol{\nabla}^2\bold{u}(t)=\bold{b}(t),
\end{equation}
where $\bold{m}$ is the square of slowness, $\bold{u}(t)$ is the discretized wavefield and $\ddot{\bold{u}}(t)$ denotes its second time derivative, $\boldsymbol{\nabla}^2$ is the Laplacian operator with respect to the spatial coordinates, and $\bold{b}(t)$ is the source term.
Then, minimization of the source residual, equation \ref{main_PMS2}, leads to
\begin{equation} \label{B-2}
\bold{m}_{k+1}=\frac{\int_0^{\tau} \ddot{\bold{u}}^e_k(t)\circ[\nabla^2\bold{u}^e_k(t)+\bold{b}_k^-(t)]dt}{\int_0^{\tau} |\ddot{\bold{u}}^e_k(t)|^2 dt},
\end{equation}
where ${\bold{u}}^e_k(t)$ obeys the following equation (see equation \ref{main_PMS2} for $\bold{b}_k^-(t)$ and $\bold{b}_k^+(t)$)
\begin{equation} \label{PDE_uk}
\bold{m}_k\circ \ddot{\bold{u}}^e_k(t)-\boldsymbol{\nabla}^2\bold{u}^e_k(t)=\bold{b}_k^+(t).
\end{equation}
Using the expression for $\boldsymbol{\nabla}^2\bold{u}^e_k(t)$, equation \ref{B-2} can also be written as
\begin{align} \label{B-3}
\bold{m}_{k+1} &=\frac{\int_0^{\tau} \ddot{\bold{u}}^e_k(t)\circ[\bold{m}_k\circ \ddot{\bold{u}}^e_k(t)-\alpha_k\bold{v}^e_k(t)]dt}{\int_0^{\tau} |\ddot{\bold{u}}^e_k(t)|^2 dt},\\
&= \bold{m}_k + \alpha_k\delta\bold{m}_{k}, \label{newton}
\end{align}
where
\begin{align} \label{EFWI_adjoint_state}
 \delta\bold{m}_{k} = -\frac{\int_0^{\tau} \ddot{\bold{u}}^e_k(t)\circ\bold{v}^e_k(t)dt}{\int_0^{\tau} |\ddot{\bold{u}}^e_k(t)|^2 dt},
\end{align}
in which $\bold{v}^e_k(t)$ is the backward-wavefield obtained by back-propagation of the adjoint source $\bold{y}_{k}+ \delta\bold{d}^r_k$.
As seen, the new algorithm simply adds back the residual to improve the estimate.

\subsubsection{Computational Cost}
The computational burden of DRI algorithm lies in the simulation of propagating wavefields using a finite difference method. More precisely, two fields ${\bold{u}}^e_k$ and ${\bold{v}}^e_k$ must be simulated and need to be available at the same time level to compute the integral in equation \ref{EFWI_adjoint_state}. 
One can compute these integrations by first computing the forward wavefield  ${\bold{u}}^r_k$ (solution of one forward problem) and storing it in the computer memory. This first step also delivers the data residual $\delta\bold{d}^r_k$.
Then the residual field $\delta{\bold{u}}_k$ is computed by back-propagation of the data residual  followed by forward propagation of the resulting backward field (solution of one forward plus one backward problems), during which the step length $\alpha_k$, equation \ref{alpha}, and the wavefield ${\bold{u}}^e_k={\bold{u}}^r_k+\alpha_k\delta{\bold{u}}_k$ are computed. Finally, the backward field $\bold{v}^e_k$ is computed by back-propagating the adjoint source  
$\bold{y}_{k}+\delta\bold{d}^r_k$ (solution of one backward problem) and the imaging condition in equation \ref{EFWI_adjoint_state} is computed simultaneously.
 Therefore, this approach requires one to solve 2 + 2 wave equations at each iteration. 
In practice, storing the wavefield is memory intensive. To mitigate this issue, one can use the compression and reconstruction methods to decrease demanded memory \citep[e.g., ][]{Boehm_2015_WCA,Gholami_2019_EWR,Zand_2019_CIR}. An alternative method, which is less memory demanding but adds the cost of extra wave equation solves, is to store only the spatial boundary conditions and final time slice and then recomputing the wavefields \citep{Gauthier_1986_TDN,Symes_2007_RTM}.

\section{NUMERICAL EXAMPLES}

\subsection{Camembert model test}
\begin{figure}
\center
\includegraphics[width=1\columnwidth]{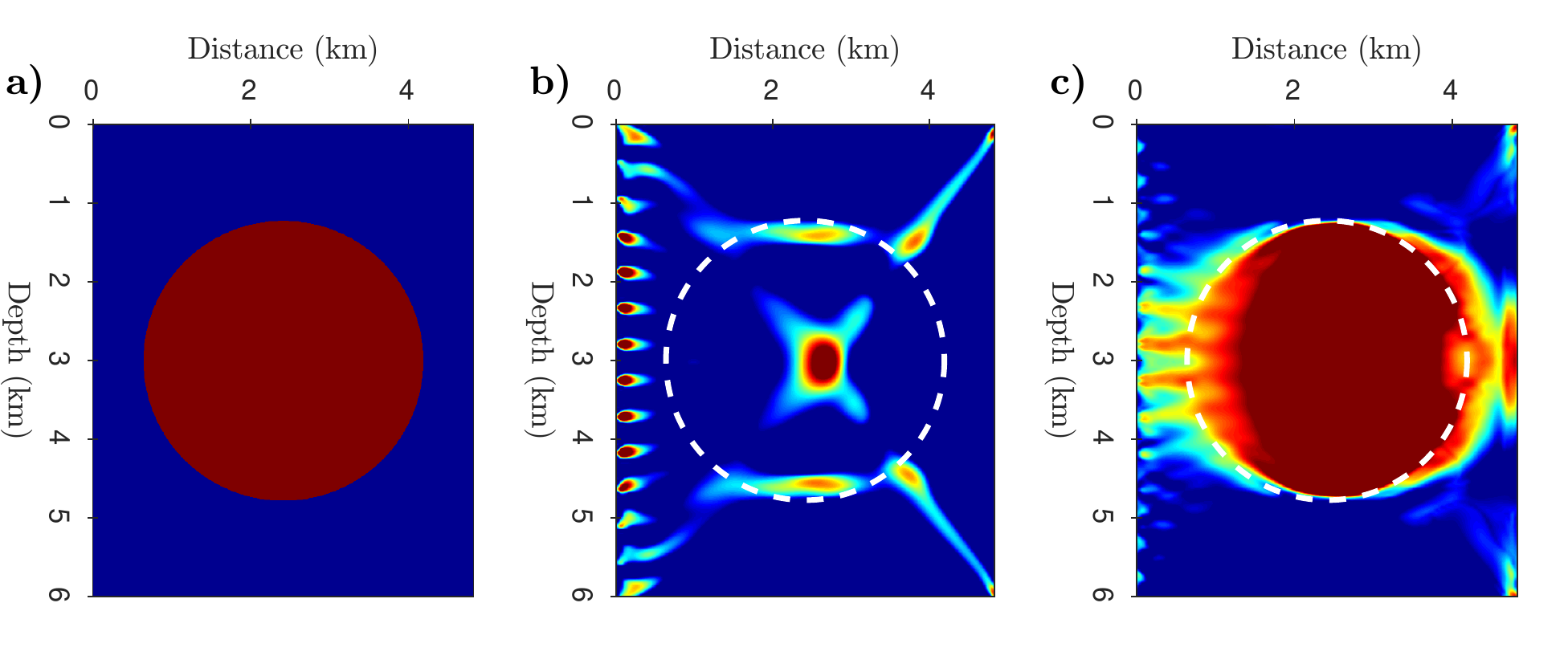}
\caption{Camembert model example. (a) True velocity, and inverted velocity models using (b) FWI and (c) DRI.}
\label{Camembert}
\end{figure}

\begin{figure}
\center
\includegraphics[width=0.5\columnwidth]{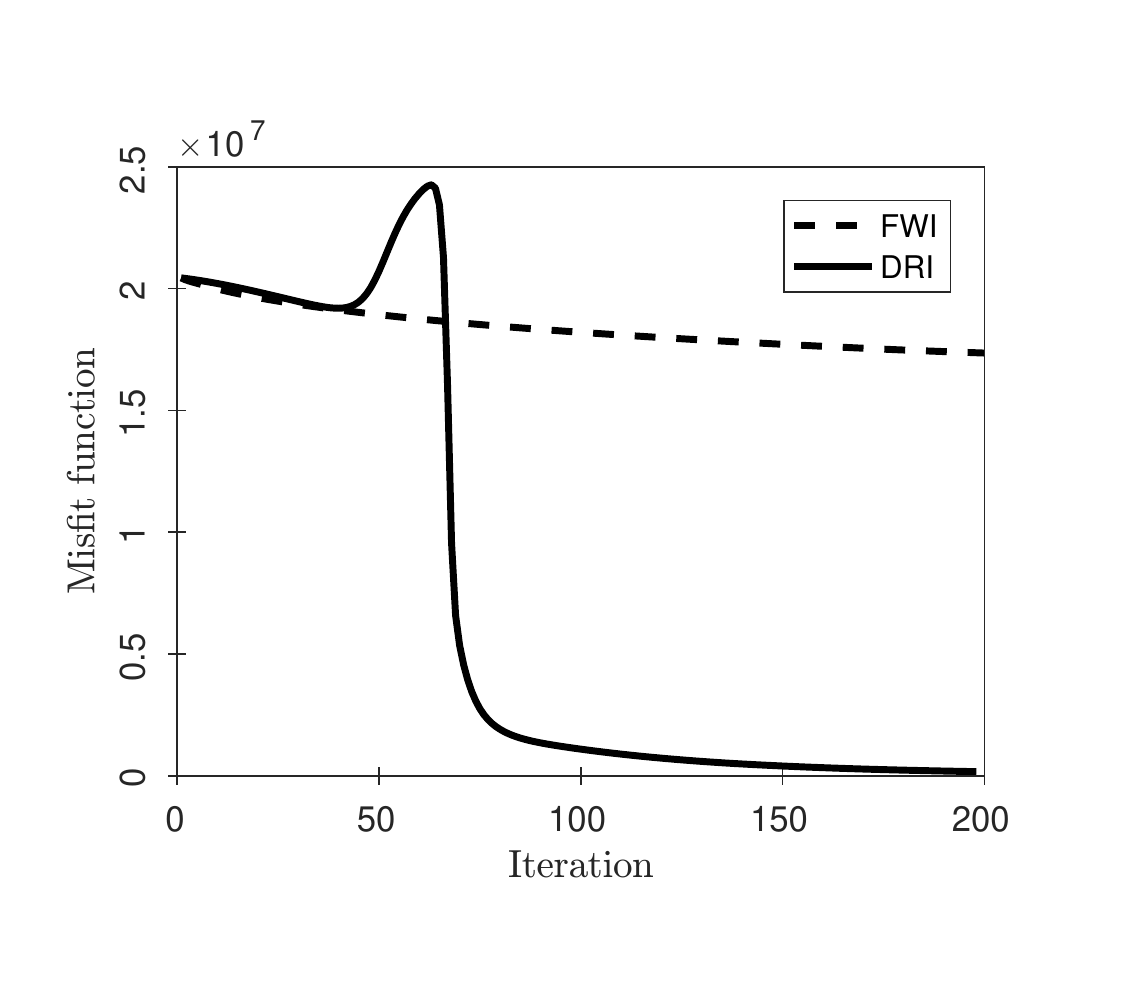}
\caption{Trajectory of the misfit function for the Camembert model example in Figure \ref{Camembert}.}
\label{Camembert_convergence}
\end{figure}

In order to show the robustness of the  proposed method with respect to the initial model, we consider the simple while challenging ``Camembert" model \citep{Gauthier_1986_TDN}. 
The subsurface model contains a circular anomaly of velocity 4.6 km/s embedded in a homogeneous background of velocity 4.0 km/s (Figure \ref{Camembert}a). 
The dimensions of the model are 4.8 km in distance and 6 km in depth, and the grid spacing is 35.5 m.
The crosshole acquisition consists of 14 equally spaced sources of 10 Hz Ricker wavelet in the left side of the model and 170 equally spaced receivers deployed vertically on the opposite side.
This model was previously used by \citet{Engquist_2020_OTB} for an application of FWI based upon
an optimal-transport distance.
Due to the large diameter of the anomaly the inverse problem is nonlinear and thus the initial model plays an important role in the convergence of the local optimization algorithms.   
We start the inversion from the homogeneous background model to meet the nonlinear regime of the FWI and perform one cycle of inversion with the 10 Hz Ricker wavelet.
The results obtained using FWI and the DRI are presented in Figures \ref{Camembert}c and \ref{Camembert}d, respectively. Figure  \ref{Camembert_convergence} shows the trajectory of the misfit function versus iteration for both methods.
The result of FWI is quite different from the true one, which indicates that the method suffers from severe cycle skipping due
to the rough initial model and 10 Hz frequency inversion. In contrast, the DRI mitigates efficiently the cycle skipping. This is clearly seen from the reconstructed model which is quite close to the ground truth.

\subsection{Checkerboard model test}
\begin{figure}
\center
\includegraphics[width=1\columnwidth]{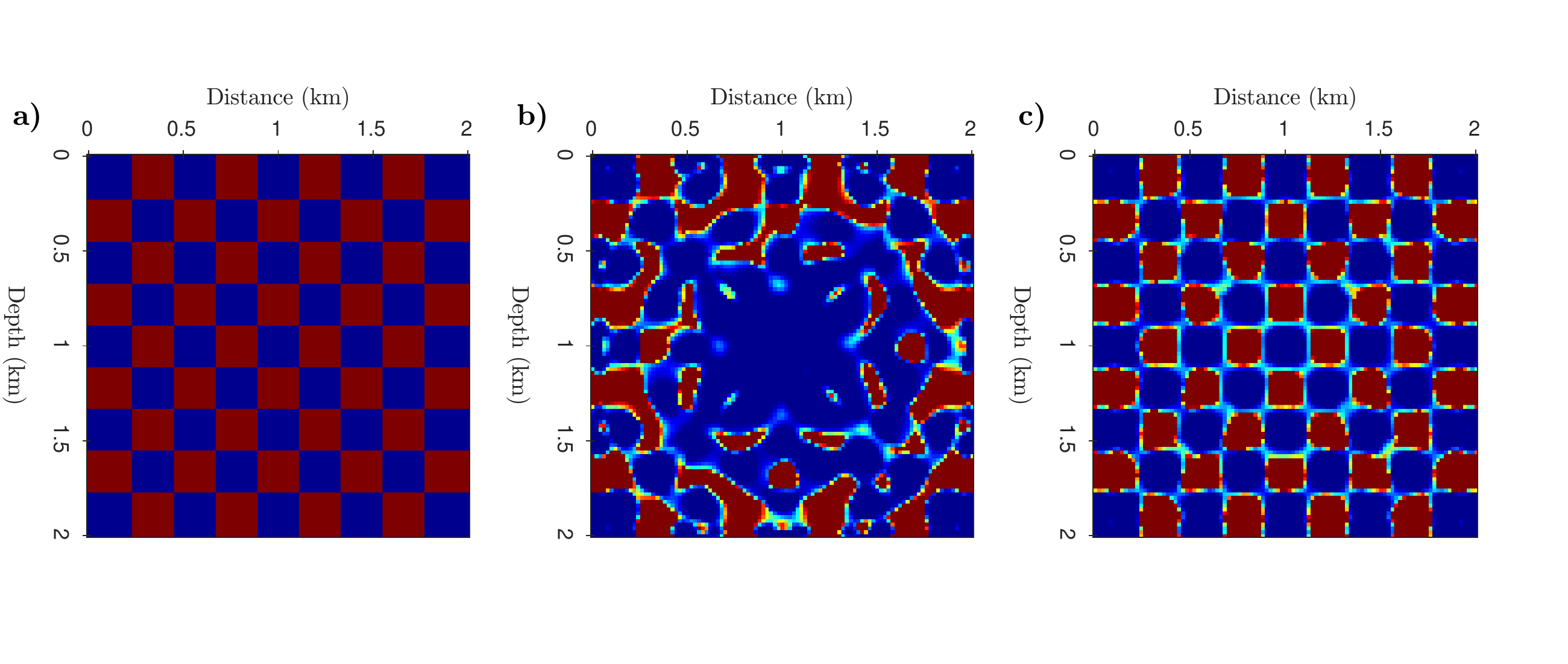}
\caption{A checkerboard model example. (a) True velocity, and inverted velocity models using (b) FWI and (c) DRI.}
\label{Checkerboard}
\end{figure}

\begin{figure}
\center
\includegraphics[width=0.5\columnwidth]{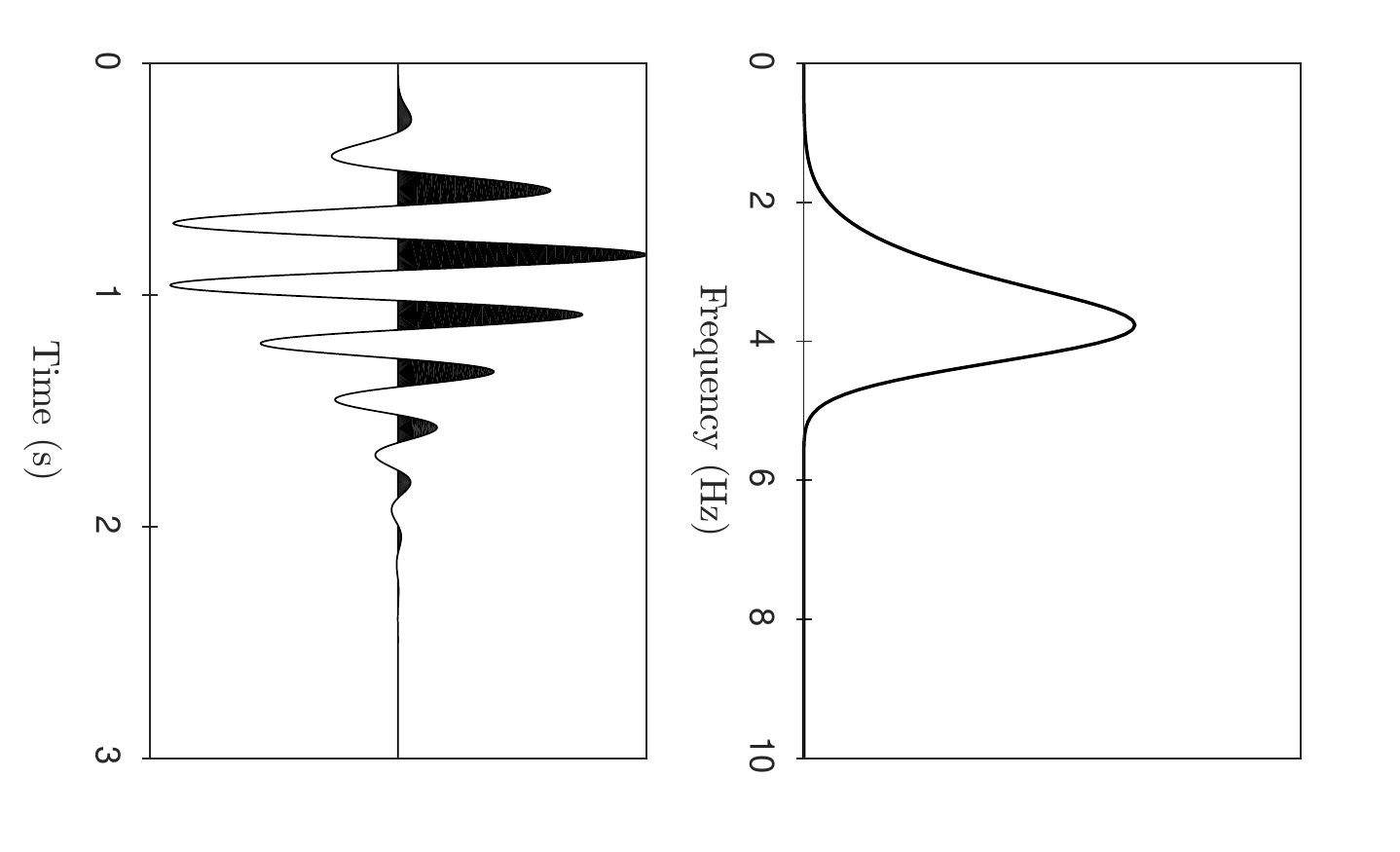}
\caption{Source wavelet in the time and frequency domains used for the checkerboard example.}
\label{Checkerboard_wavelet}
\end{figure}

\begin{figure}
\center
\includegraphics[width=0.7\columnwidth]{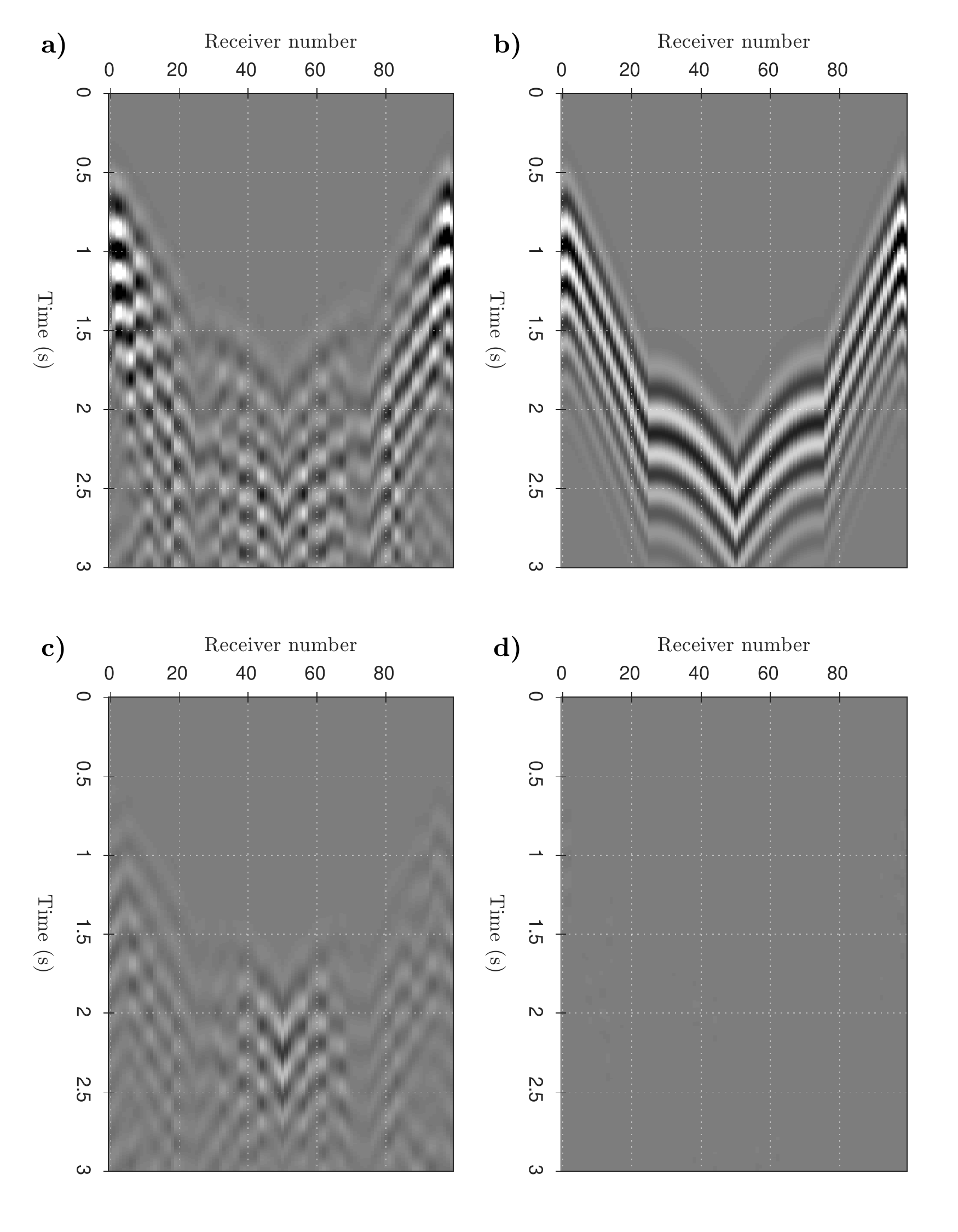}
\caption{A data-set of the checkerboard model. Data-set computed in the exact model (a) and in the initial homogeneous model (b). Residual data at the convergence point for (c) the FWI and (d) the DRI methods. }
\label{Checkerboard_data}
\end{figure}

\begin{figure}
\center
\includegraphics[width=0.5\columnwidth]{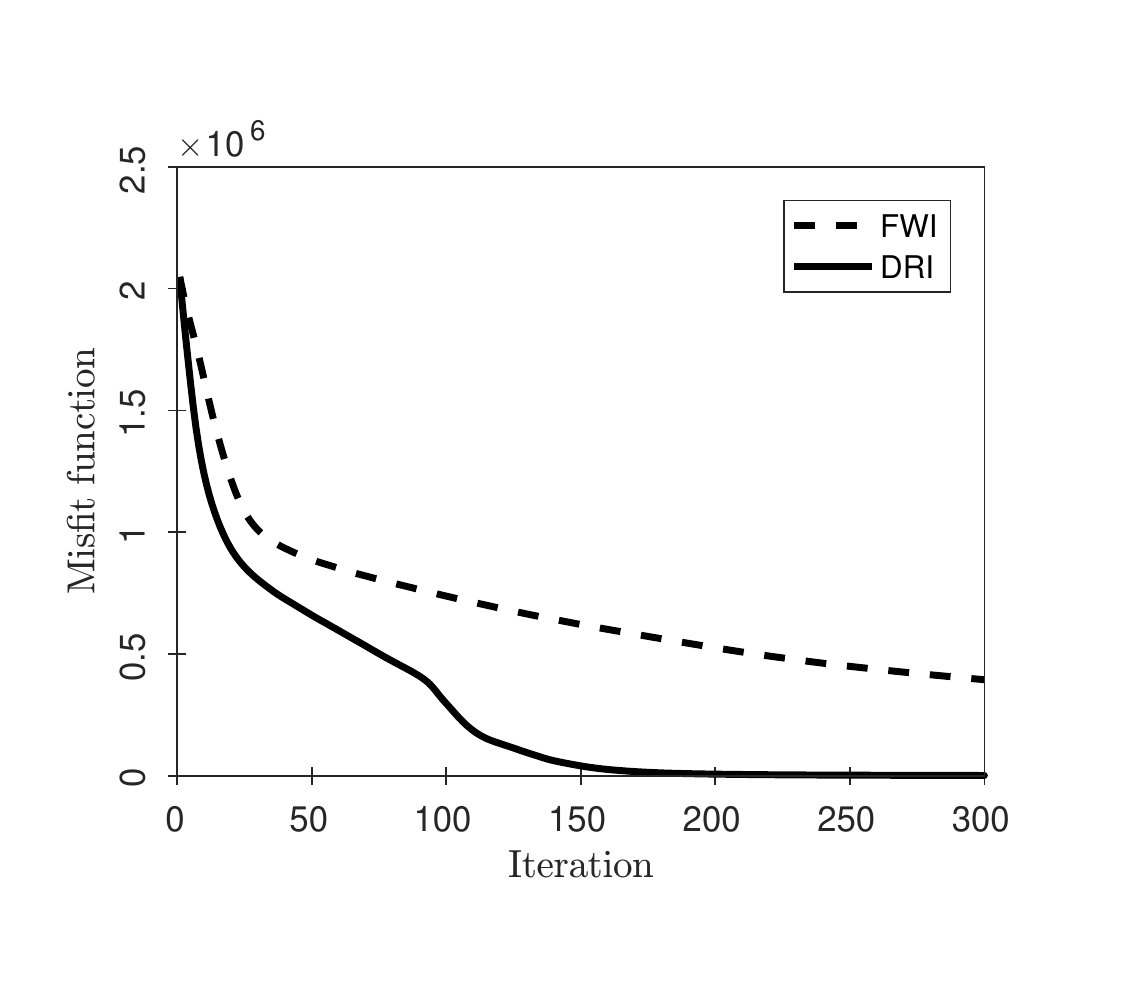}
\caption{Trajectory of the misfit function for the checkerboard model example in Figure \ref{Checkerboard}.}
\label{Checkerboard_convergence}
\end{figure}

We design a checkerboard model composed of a homogeneous background of velocity 1.5 km/s and checkerboard
velocity perturbations of velocity 4.0 km/s and dimension 200 m by 200 m (Figure \ref{Checkerboard}a).
We use a full acquisition with 92 equally-spaced sources placed all around the domain with 50 m distance from the model boundary. 
Each source is associated with a total of 100 equally spaced receivers placed on the model boundary.
The model is discretized over a 101 by 101 grid with a spatial step of 20 m.
The source wavelet is a band-pass filtered Ricker wavelet between 2.5 Hz and 5 Hz with a peak frequency of 3.75 Hz (Figure \ref{Checkerboard_wavelet}). 
The initial model is the background velocity.
Figures \ref{Checkerboard_data}a and \ref{Checkerboard_data}b show two data-sets computed with the exact velocity model and with the initial homogeneous-velocity model.
We compare the inversion results obtained performing 300 iterations (Figure \ref{Checkerboard_convergence}) for both FWI and DRI and the corresponding model estimates are shown in Figures \ref{Checkerboard}b and \ref{Checkerboard}c.
The final data residuals corresponding to both methods are also shown in Figures \ref{Checkerboard_data}c and \ref{Checkerboard_data}d.
For this test, again, unlike FWI, the DRI reconstructed a model which is quite close to the ground truth.

\subsection{Marmousi II Model test}
\begin{figure}
\center
\includegraphics[width=0.5\columnwidth]{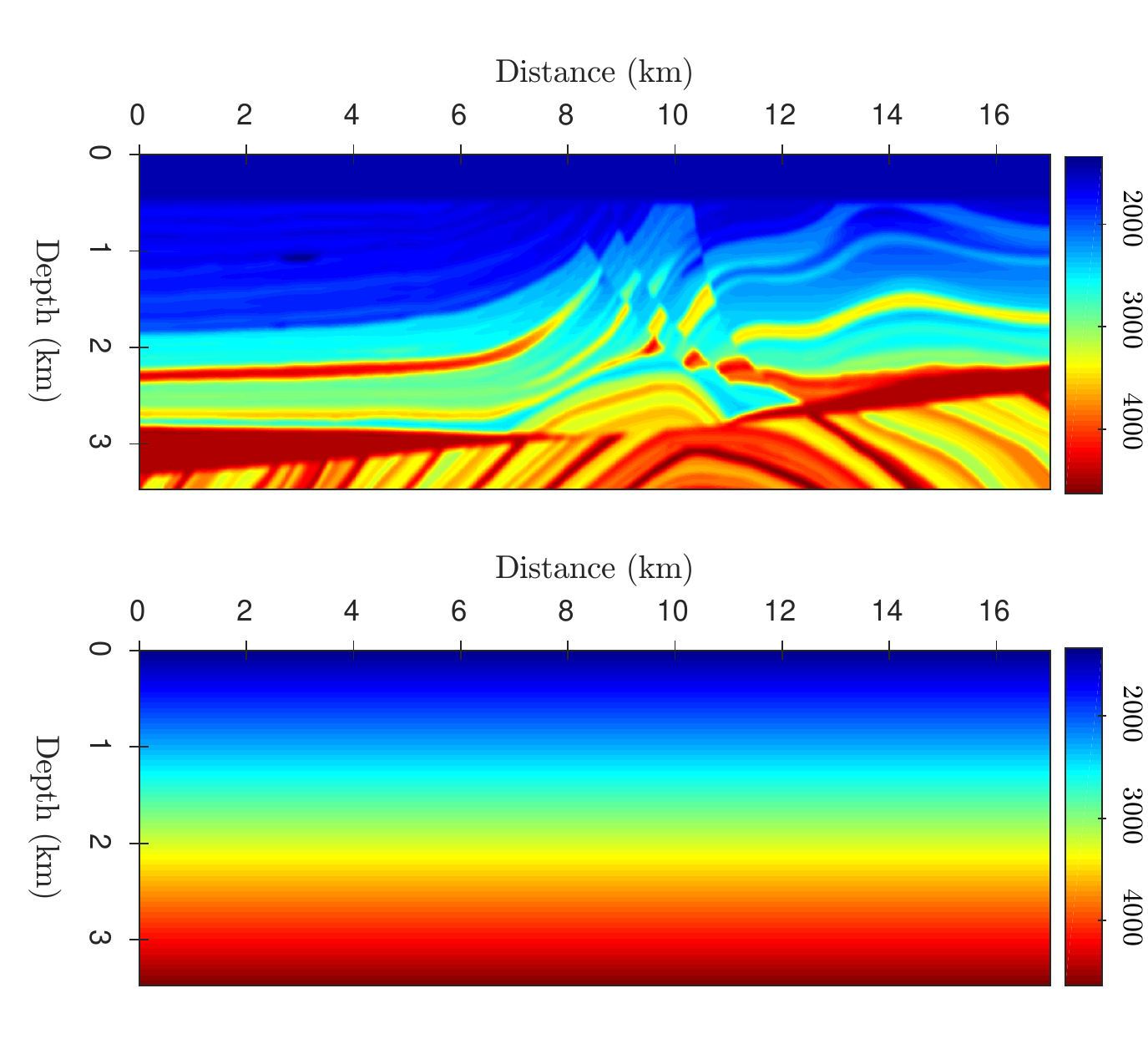}
\caption{Marmousi II velocity model (top) and 1D initial model (bottom).}
\label{Marmousi_vel}
\end{figure}

\begin{figure}
\center
\includegraphics[width=0.5\columnwidth,trim={0 5cm 0 0cm},clip]{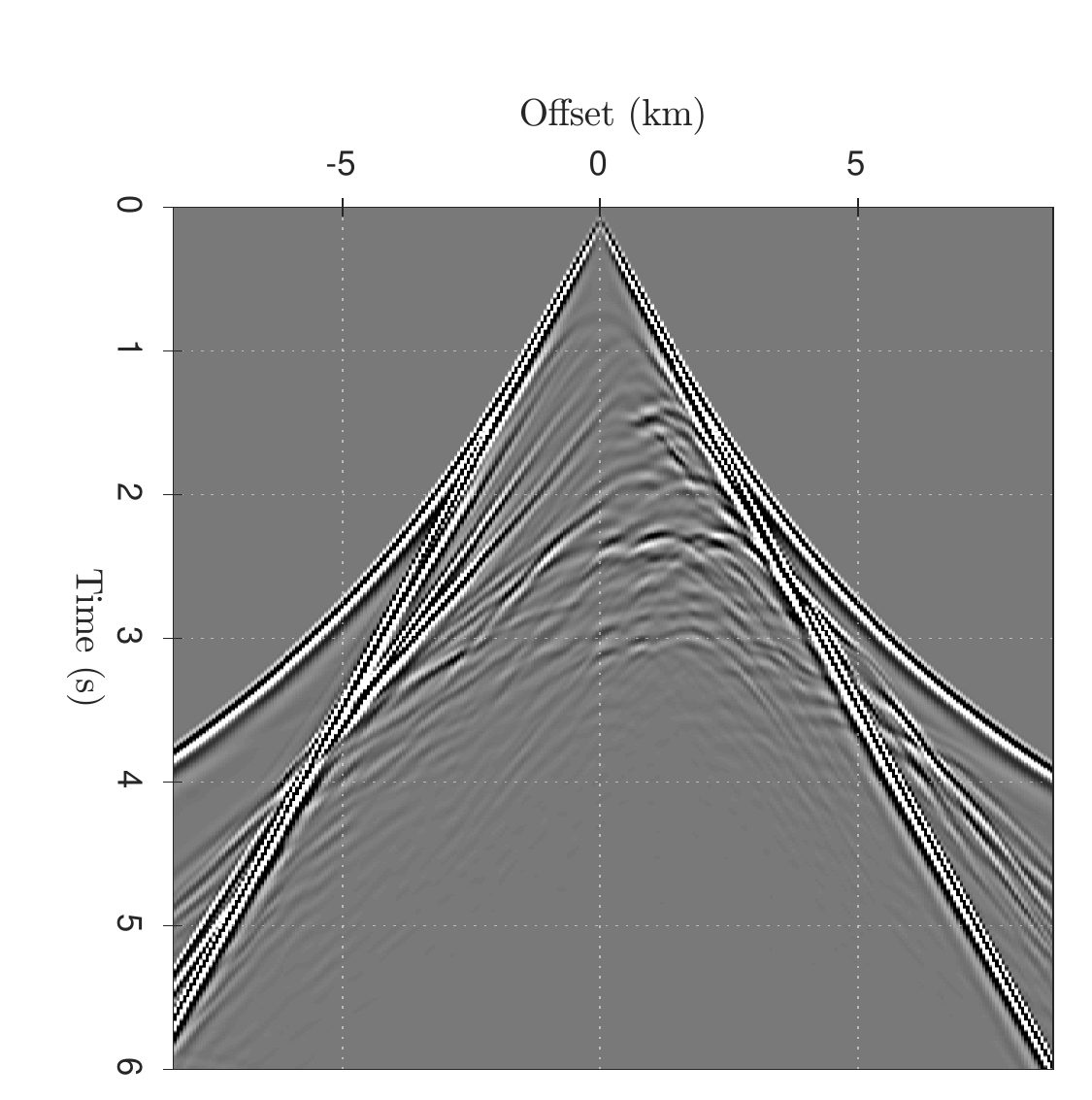}
\caption{Residual data of Marmousi model for a source at distance 8250 m.}
\label{Marmousi_residual}
\end{figure}

\begin{figure}
\center
\includegraphics[width=0.5\columnwidth,trim={0 0cm 0 0cm},clip]{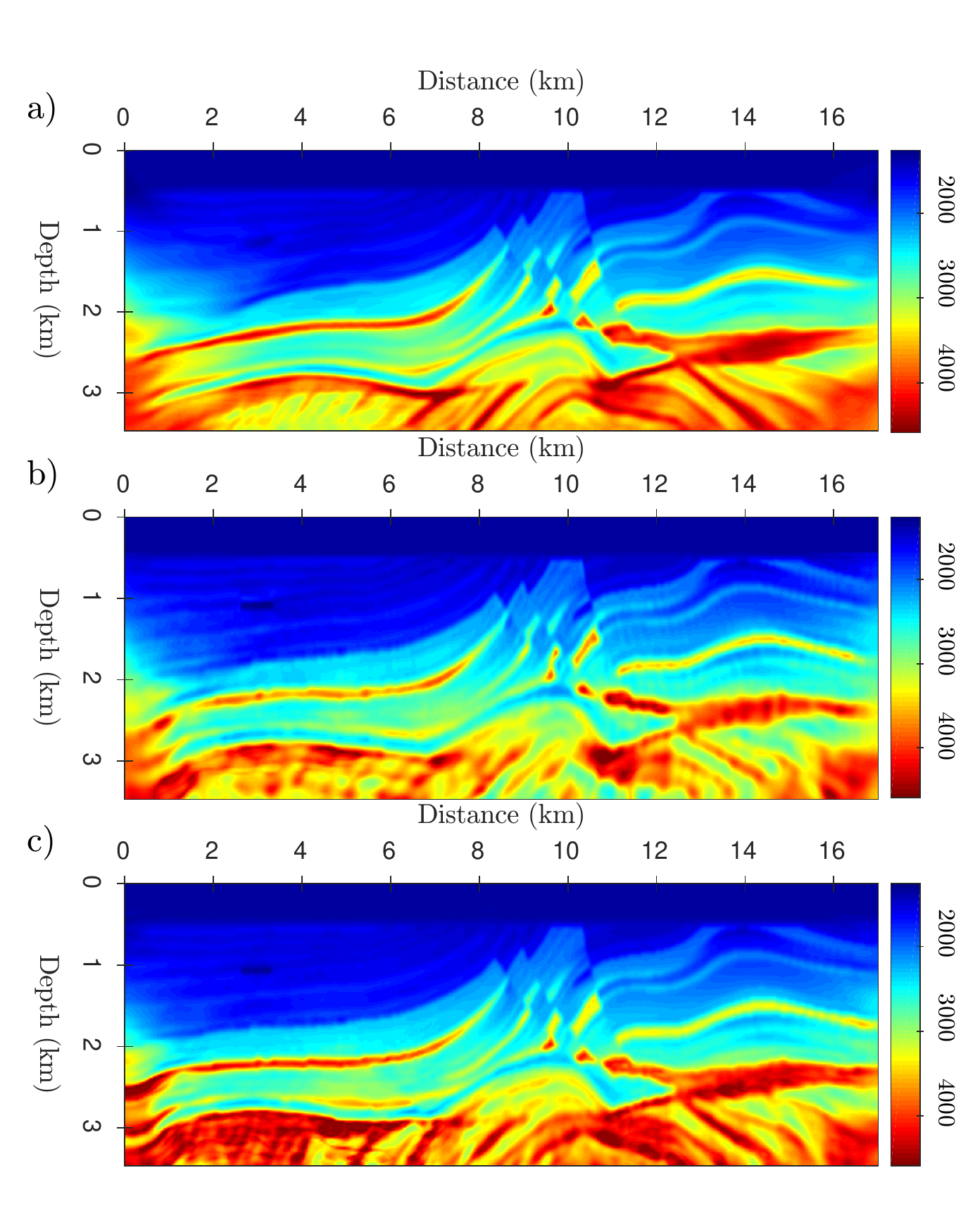}
\caption{Marmousi II velocity model obtained by (a) time-domain DRI, (b) frequency-domain DRI, (c) frequency-domain IR-WRI.}
\label{Marmousi_vel_estimate}
\end{figure}

We now consider a more realistic multisource experiment using the Marmousi II velocity model (Figure \ref{Marmousi_vel}). 
The model is 17 km long and 3.5 km deep. The fixed spread acquisition consists of equally-spaced 91 sources and equally-spaced 273 receivers on top side of the model. 
We build reference data using a highpass filter minimum-phase Ricker wavelet with no signal below 2.5 Hz. The total recording time is 7.5 seconds with a sample internal 3 ms.
The initial model is  a 1D velocity model linearly increasing in depth from 1.5 km/s to 4.5 km/s.
Figure \ref{Marmousi_residual} shows the residual gather for the source at distance 8.5 km. 
Obviously classical FWI will stuck in a local optimum due to severe cycle skipping caused by the 1D initial model and the absence of low-frequencies in the data \citep[see][ Their Figure 6]{Aghamiry_2019_IWR} thus we don't show the corresponding results here. 
The inversion result of the DRI is shown in Figure \ref{Marmousi_vel_estimate} which clearly shows that the method mitigated efficiently the cycle skipping. The reconstruction is quite close to the true model in most parts of the model.
In order to compare the DRI algorithm with the IR-WRI algorithm we performed both algorithm in the frequency domain to invert the data corresponding to the Marmousi II. The results obtained by both algorithms are shown in Figures \ref{Marmousi_vel_estimate}b-\ref{Marmousi_vel_estimate}c and as seen the estimated velocity models are similar which confirm the accuracy of the approximation we considered in DRI algorithm, given in equation \ref{hess_aprox}.

\subsection{2004 BP Model test}
\begin{figure}
\center
\includegraphics[width=0.5\columnwidth,trim={0 0cm 0 0cm},clip]{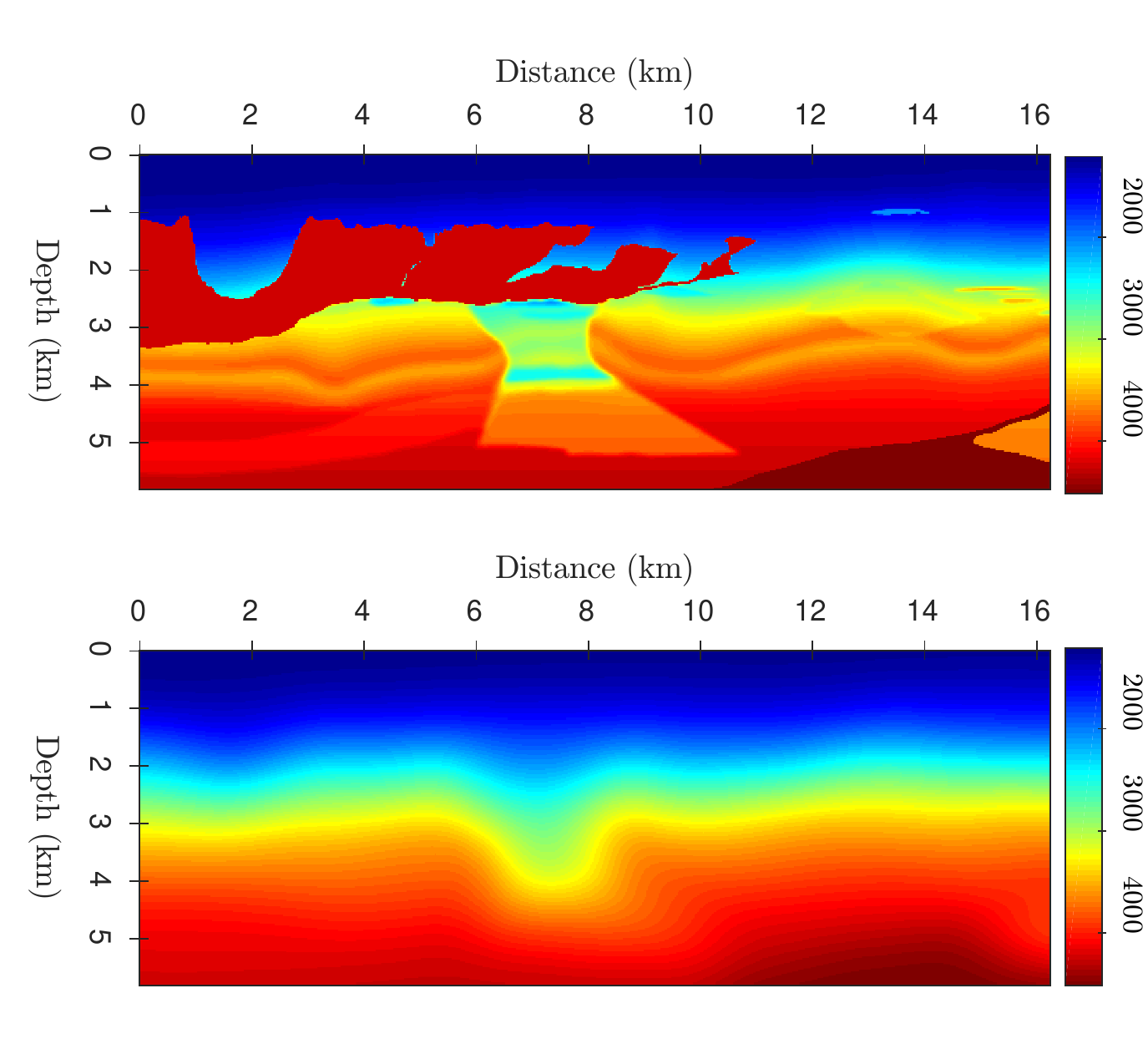}
\caption{2004 BP velocity model (top) and initial model (bottom).}
\label{BP_vel}
\end{figure}

\begin{figure}
\center
\includegraphics[width=1\columnwidth,trim={0 0cm 0 0cm},clip]{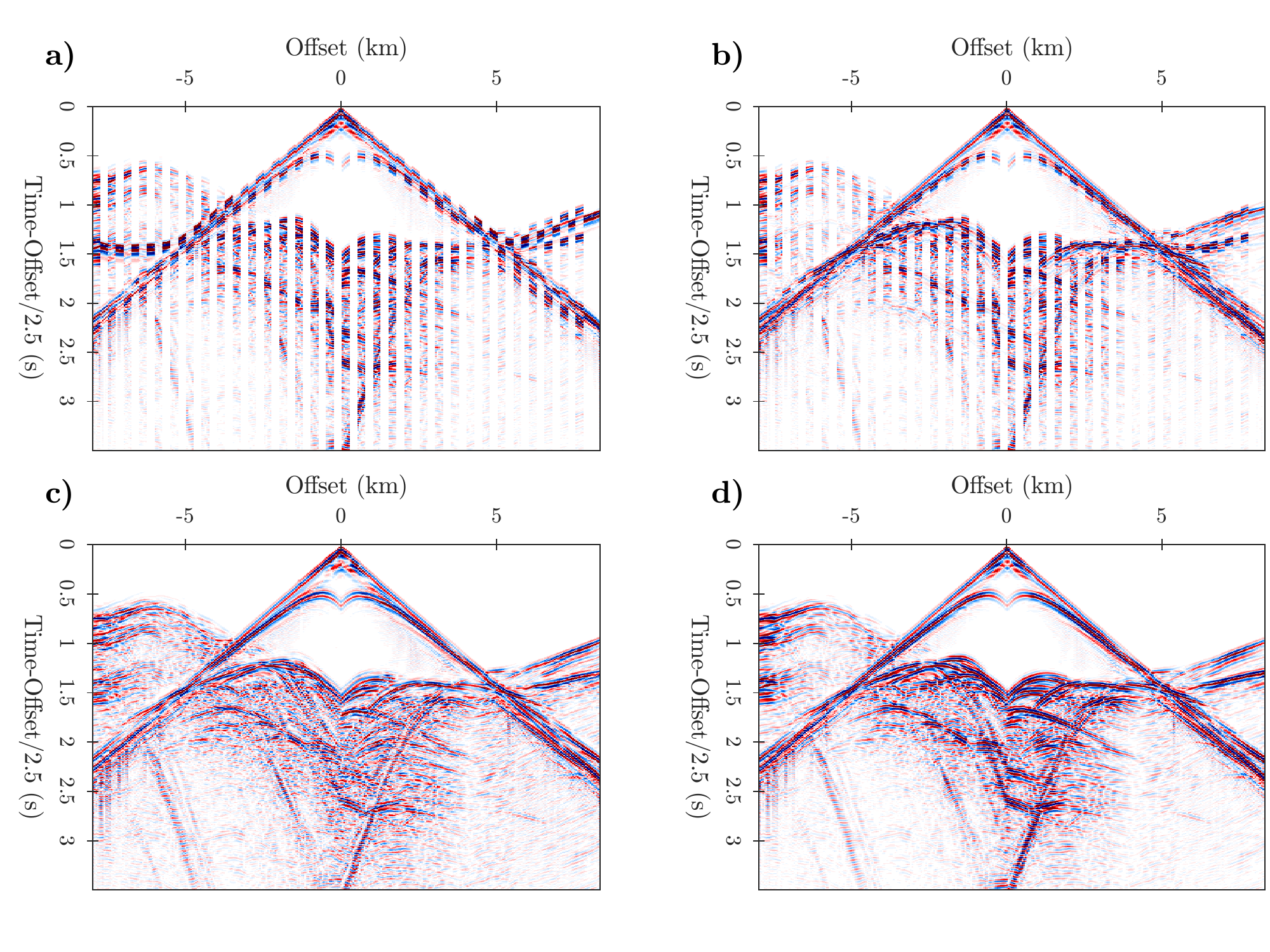}
\caption{Comparison between a common-shot gather computed in (a) the initial model, (b) the estimated model at first iteration, and (c) the estimated model at final iteration. (d) The gather computed in the true 2004 BP velocity model.}
\label{BP_data}
\end{figure}

\begin{figure}
\center
\includegraphics[width=0.5\columnwidth,trim={0 0cm 0 6cm},clip]{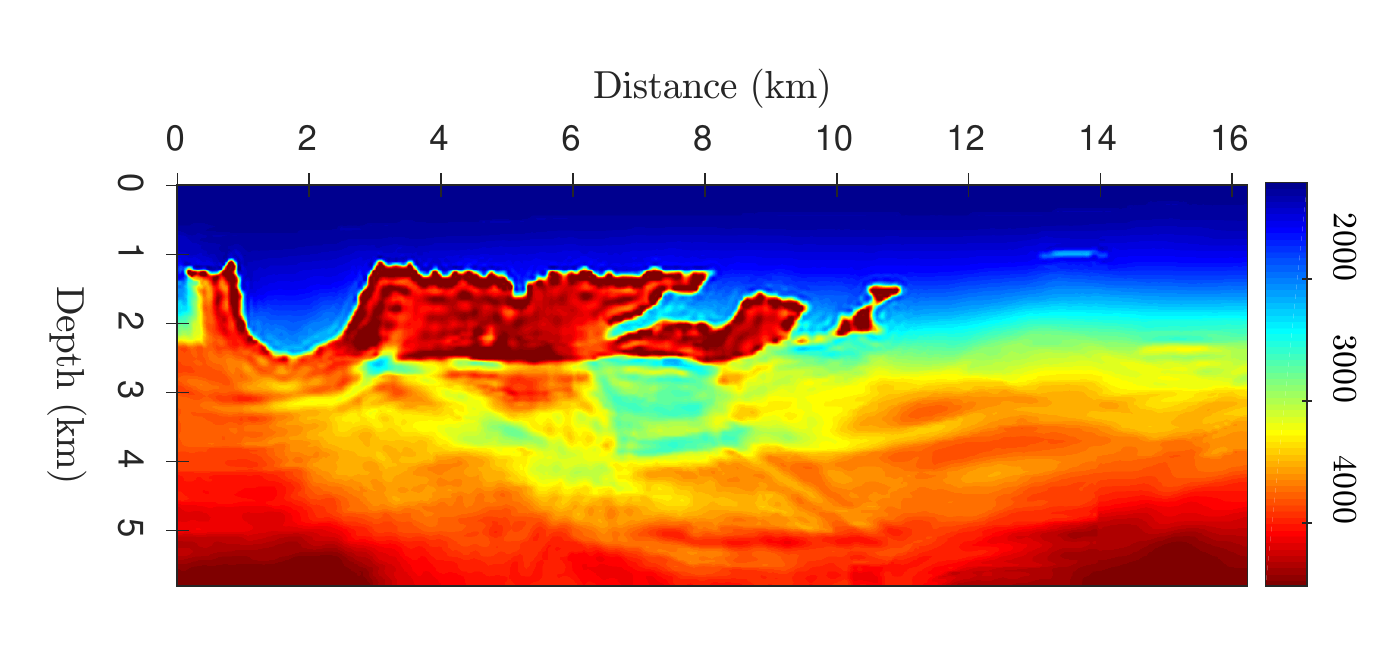}
\caption{2004 BP  velocity model obtained by the proposed algorithm.}
\label{BP_vel_estimate}
\end{figure}

Finally, we assess the DRI method against a target of the challenging 2004 BP salt model  when a crude starting model is used (Figure \ref{BP_vel}).
The selected target corresponds to the left part of the 2004 BP salt model and was previously used by \citet{Metivier_2015_MTM} for an application of FWI based upon an optimal-transport distance and by \citet{Aghamiry_2019_IWR} for an application of frequency-domain WRI based upon the augmented Lagrangian formulation (IR-WRI).
The subsurface model is 16.25 km wide and 5.825 km deep. 
The fixed spread acquisition consists of equally-spaced 109 sources and equally-spaced 326 receivers on top side of the model. 
We build reference data using a highpass filter minimum-phase Ricker wavelet with no signal below 2.5 Hz. The total recording time is 7.5 seconds with a sample internal 2 ms.
We used a smoothed version of the true velocity model without the salt as initial model (Figure \ref{BP_vel}).
A comparison between a common-shot gather computed in the true and initial models shows severe traveltime mismatches (Figure \ref{BP_data}) that drive the classical FWI to a local minimum \citep[see][ Their Figure 12c]{Aghamiry_2019_IWR}.
The DRI, however, mitigated efficiently the cycle skipping as it can be seen from the reconstructed model shown in Figure \ref{BP_vel_estimate}.
As seen from this figure, the precise geometry of the salt body and the subsalt low-velocity anomalies are mostly recovered with a quite high resolution.

\section{Conclusions}
Time domain full waveform inversion (FWI) is appealing since the wave equation can be solved efficiently by explicit time-stepping but this method require sufficiently accurate initial models which may not be available in practice. 
Extended FWI with wavefield reconstruction is robust with respect to the initial model and it can be solved with a minimal sensitivity to the penalty parameter choice when using the augmented Lagrangian solver. However, its application in the time domain is a challenging task due to the absent of explicit time-stepping and because computation and memory burdens grow quickly as the size of model increases due to the inherent high-dimensional nature of the state variable.
In this paper, we proposed a new augmented Lagrangian algorithm for this task which performs in the data space, a  significantly lower dimensional subspace than the dimensionality of the original full space.
In this algorithm, the wavefield reconstruction appears as a residual data reconstruction
and the model parameters are updated by a quasi-Newton like algorithm with a search direction computed by a generalized adjoint state method. This latter is computed by zero-lag cross-correlation between the data-assimilated wavefields and adjoint fields coming from back-propagation of all previous data residuals added back to the current residual vector.   
Approximate update of the residual reconstruction by a gradient descent method resulted to an efficient algorithm which only requires four wave equation solves (two forward solves plus two backward solves) at each iteration. 
Different numerical tests including Marmousi II and 2004 BP velocity models confirm the efficiency of this algorithm for time-domain applications similar to the classical FWI while benefiting from the excellent robustness of the extended FWI.

\section{Acknowledgments}  
This study was partially funded by the WIND consortium (\textit{https://www.geoazur.fr/WIND}), sponsored by Chevron, Shell, and Total. This study was granted access to the HPC resources of SIGAMM infrastructure (http://crimson.oca.eu), hosted by Observatoire de la C\^ote d'Azur and which is supported by the Provence-Alpes C\^ote d'Azur region, and the HPC resources of CINES/IDRIS/TGCC under the allocation A0050410596 made by GENCI."

 \appendix
\section{\\Appendix A}
\label{appendixa}
\begin{figure}
\center
\includegraphics[width=1\columnwidth]{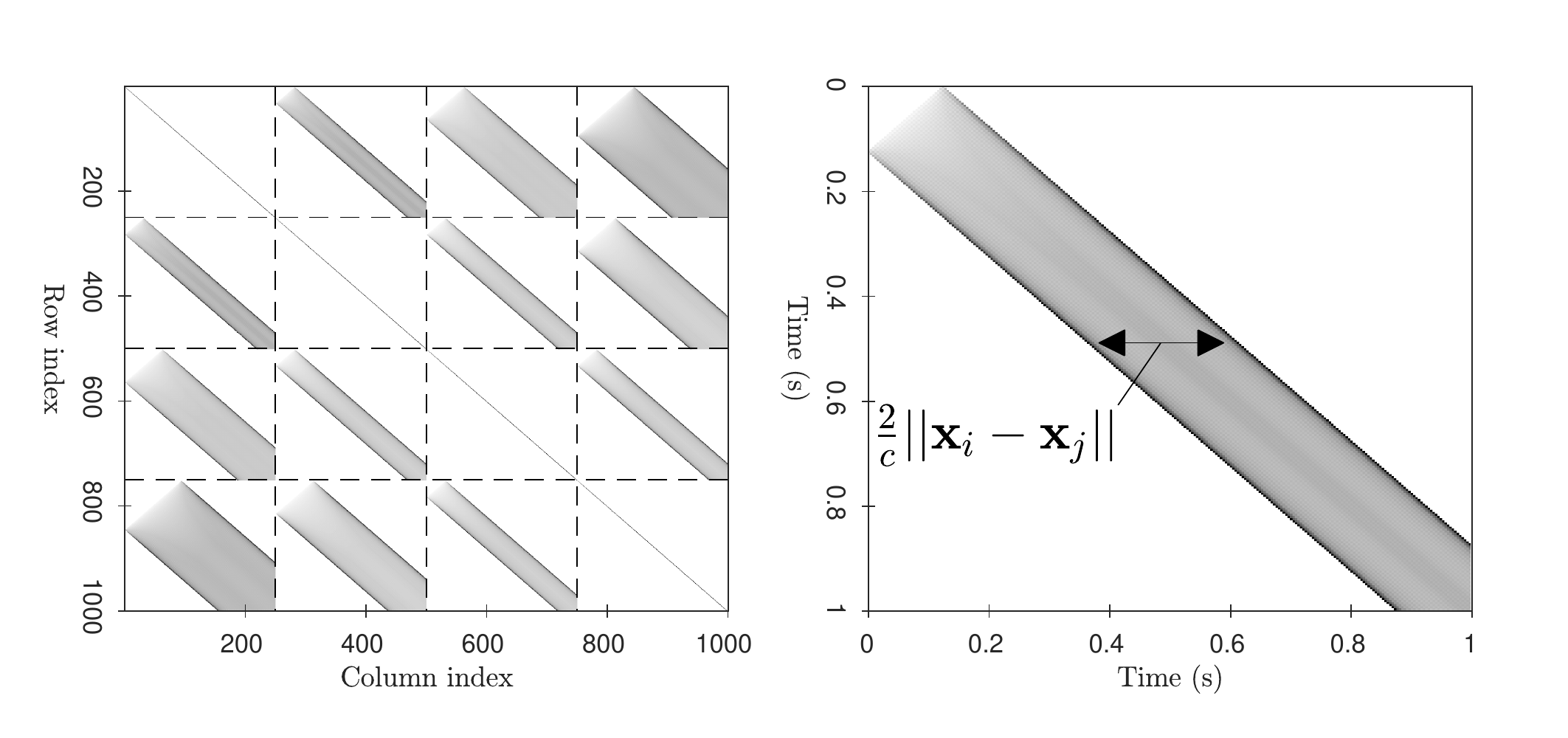}
\caption{(Left) Q matrix computed for a homogeneous model of velocity 2000 m/s and for 5 receivers each of 251 time samples. The distance between receivers is 200 m. (Right) block (1,2) of matrix Q where the dashed lines show the theoretical band.}
\label{Qmat}
\end{figure}

Consider the acoustic wave equation in a homogeneous 3D medium
\begin{equation} \label{PDE_homo}
\frac{1}{c^2}\frac{\partial^2}{\partial t^2} u(t,\bold{x})-\nabla^2u(t,\bold{x})=b(t,\bold{x}),
\end{equation}
in which the acoustic velocity $c$ is a constant.
The solution of equation \ref{PDE_homo} can be written in the following integral form:
\begin{equation} \label{green}
u(t,\bold{x}) = \iint g(t-t',\bold{x}-\bold{x}') b(t',\bold{x}') dt' d\bold{x}',
\end{equation}
where 
\begin{equation}
g(t,\bold{x}) = \frac{\delta(t - \|\bold{x}\|_2/c)}{\|\bold{x}\|_2},
\end{equation}
are the Green's functions and $\delta$ denotes the delta function. Writing equation \ref{green} in an operator notation as $u = A^{-1}b$ then the goal is to build 
\begin{equation} \label{Q}
Q = PA^{-1}A^{-T}P^T.
\end{equation}
where $P$ is the restriction matrix at the receiver locations. 
In order to build $Q$ we need to compute the adjoint wavefield for each column of $P^T$ followed by computing the forward wavefield while using the ajoint field as the source and then sampling the result at receiver location. 
Each column of $P^T$ can be written as 
\begin{equation} \label{delta}
\delta(t-\tau,\bold{x}-\bold{x}_j)
\end{equation}
in which $\bold{x}_j$ denotes the location of the $j$th receiver and $\tau$ denotes the time lag.
The adjoint wavefield for the delta function in equation \ref{delta} is 
\begin{align} \label{adjgreen}
v(t,\bold{x}) &= \iint g(t'-t,\bold{x}'-\bold{x}) \delta(t'-\tau,\bold{x}'-\bold{x}_j) dt' d\bold{x}'\\
&=g(\tau-t,\bold{x}_j-\bold{x}).
\end{align}
Then applying the forward operator on this adjoint wavefield gives
\begin{align} 
u(t,\bold{x}) &= \iint g(t-t',\bold{x}-\bold{x}') v(t',\bold{x}') dt' d\bold{x}'\\
&= \iint g(t-t',\bold{x}-\bold{x}') g(\tau-t',\bold{x}_j-\bold{x}') dt' d\bold{x}'.
\end{align}
Finally, sampling this wavefield at the $i$th received location, thus counting for the matrix $P$, gives
\begin{align} 
Q(t,\tau,\bold{x}_i,\bold{x}_j) &= \iint g(t-t',\bold{x}_i-\bold{x}') g(\tau-t',\bold{x}_j-\bold{x}') dt' d\bold{x}'.
\end{align}
Substituting for the Green's functions from equation \ref{green} and performing the integration over $t'$ gives
\begin{align} \label{Qint}
Q(t,\tau,\bold{x}_i,\bold{x}_j) &=\int \frac{\delta(t-\tau+ \frac{1}{c}[\|\bold{x}_j-\bold{x}'\|_2 - \|\bold{x}_i-\bold{x}'\|_2])}{\|\bold{x}_i-\bold{x}'\|_2\|\bold{x}_j-\bold{x}'\|_2}
 d\bold{x}',
\end{align}
In order to compute this integral, we note that for each receiver pair $i$ and $j$ and a point $\bold{x}'$ the $Q$ matrix includes a linear event parallel to the main diagonal described by 
\begin{equation}
t = \tau + \frac{1}{c}[\|\bold{x}_j-\bold{x}'\|_2 - \|\bold{x}_i-\bold{x}'\|_2],
\end{equation}
in which 
the second term at right hand side
determines the amount of shift from the main diagonal. 
%
From the backward triangle inequality we have that
\begin{equation}
\|\bold{x}_j-\bold{x}'\|_2 - \|\bold{x}_i-\bold{x}'\|_2 \leq \|\bold{x}_i-\bold{x}_j\|_2,
\end{equation}
which implies that the maximum shift is $\frac{1}{c} \|\bold{x}_i-\bold{x}_j\|_2$.
Summing over all points $\bold{x}'$ thus leads to a symmetric banded matrix with a band centred on the main diagonal. Figure \ref{Qmat} shows the Q matrix for  a homogeneous model of velocity 2000 m/s for 5 equally spaced receivers each of 251 time samples. The distance between receivers is 200 m. The band diagonal structure of each block is clearly seen and the bandwidth follows the theoretical finding. 

%

\graphicspath{{"./figures/"}}

%
%
%
%
%

\bibliographystyle{seg}

\newcommand{\SortNoop}[1]{}

\end{document}